\documentclass[12pt, reqno]{amsart}
\usepackage{paper_style}

\usepackage[csfont, paper]{std_math}

\usepackage{local}
\usepackage{Changing_Macros}

\numberwithin{equation}{section}

\makeatletter
\renewcommand{\l@subsection}{\@tocline{2}{0pt}{2.5em}{}{}}
\renewcommand{\l@subsubsection}{\@tocline{3}{0pt}{4em}{}{}}
\makeatother

\hypersetup{pdfauthor={Menny Aka and Segev Gonen Cohen}
	,pdftitle={Hecke tightness}
	,urlcolor=blue
	,citecolor=red
    ,linkcolor=blue
	,colorlinks=true
}

\begin{document}
		
\title[Equidistribution of Hecke neighbours]{Equidistribution of Hecke neighbours and \\ continued fraction patterns}
\author{Menny Aka and Segev Gonen Cohen}
\address{Department of Mathematics, ETH Zürich, Zürich, Switzerland}
\email{menny.akka@math.ethz.ch}
\email{segev.gonencohen@math.ethz.ch}

\thanks{The authors acknowledge the support of FNS (Fonds National Suisse, Equidistribution in Number Theory, no. 10.003.145).}


\date{\today.}
\begin{abstract}
We prove uniform tightness of Hecke neighbours of any probability measure of positive entropy that is invariant and ergodic under the diagonal flow. As an application, we confirm a conjecture of Einsiedler concerning a homogeneous dynamics analogue of an equidistribution theorem of Johnson and Rudolph on the torus. We also obtain applications to continued-fraction patterns in arithmetic sequences and to the \(p\)-adic Littlewood conjecture.

 Our proof studies escape of mass using dynamics on the boundary of the Chabauty compactification of the space of lattices, together with averaging operators defined by leafwise measures along horospherical subgroups.
\end{abstract}

	\maketitle


	\section{Introduction}

\noindent The starting point of this work is the following question of Einsiedler, Fishman, and Shapira. To any irrational number $\alpha$ we associate its continued fraction expansion $\alpha = [a_0;a_1,a_2,\ldots]$. We will often denote the $i^{th}$ partial quotient by $a_i(\alpha)$. Recall that $\alpha$ is \emph{badly approximable} if and only if its partial quotients are bounded; the set of badly approximable numbers has Lebesgue measure 0 but is of full Hausdorff dimension \cite{Schmidt1966}. A \emph{pattern} is a finite non-empty sequence of positive integers.

\begin{que}[{\cite[Question~5.7]{EFS2011}}]
Let $\alpha$ be badly approximable and let $w$ be a pattern. Is there an integer $n\geq1$ such that $w$ occurs infinitely often in the continued fraction expansion of $n\alpha$?
\end{que}

This question strengthens their solution of an earlier problem of Boshernitzan: the eventual bounds on the partial quotients of $n\alpha$ cannot be uniform in $n$ \cite[Theorem~1.11]{EFS2011}. This is immediate for a well approximable irrational $\alpha$, since its own partial quotients are unbounded. 

As explained in \cite[Section~5]{EFS2011}, the question is related to the distribution of the geodesic attached to $\alpha$ and its Hecke neighbours in the space of lattices. A pattern in the continued fraction expansion is coded by an open set, and multiplication of $\alpha$ by an integer corresponds to a specific choice of a Hecke neighbour, up to a shift along the geodesic. We prove two results about Hecke neighbours of measures invariant under the geodesic flow. Their applications give a positive answer to the question under an entropy assumption; in fact, it is enough to consider $\{p^k\alpha\}_{k\geq 0}$, where $p$ is any prescribed prime (see Theorem~\ref{theo:cf_patterns} below).

\subsection{Hecke neighbours: tightness and entropy}
\label{subsec:intro_hecke}

Throughout, we will take actions on the left and quotients on the right. Put
\[
 G=\SL_2(\bbR),\qquad \Gamma=\SL_2(\bbZ),\qquad X=G/\Gamma,
\]
and denote the Haar probability measure on $X$ by $m_X$. For $\varepsilon>0$, let $X(\varepsilon)$ be the compact set of unimodular lattices whose non-zero vectors all have length at least $\varepsilon$. We denote the identity element in a group by $e$. Define
\[
 a_t=\begin{pmatrix}e^t&0\\0&e^{-t}\end{pmatrix},\qquad
 u^+(s)=\begin{pmatrix}1&s\\0&1\end{pmatrix},\qquad
 u^-(s)=\begin{pmatrix}1&0\\s&1\end{pmatrix}.
\]
We write $A=\{a_t:t\in\bbR\}$ and $U^\pm=\{u^\pm(s):s\in\bbR\}$. The \emph{commensurator} of $\Gamma$ in $G$ is defined by
\begin{equation*}
    \comm_G(\Gamma):= \{r\in G: \Gamma \cap r\Gamma r^{-1}\text{ has finite index in } \Gamma \text{ and } r\Gamma r^{-1}\}.
\end{equation*}
For $r\in\comm_G(\Gamma)$ put $\Gamma_r=\Gamma\cap r\Gamma r^{-1}$ and $X_r=G/\Gamma_r$. Define the finite covering maps $p_r,q_r:X_r\to X$ by
\begin{equation*}
 p_r(g\Gamma_r)=g\Gamma,\qquad q_r(g\Gamma_r)=gr\Gamma.
\end{equation*}
Given an $A$-invariant probability measure $\mu$ on $X$, a \emph{good-$r$-lift of $\mu$} is an $A$-invariant probability measure $\lambda$ on $X_r$ with $(p_r)_*\lambda=\mu$. Set
\begin{equation*}
 \mathcal T(\mu)=\{(q_r)_*\lambda:r\in\comm_G(\Gamma),\ \lambda\text{ a good-}r\text{-lift of }\mu\}.
\end{equation*}
We refer to these measures as \emph{Hecke neighbours} of $\mu$.

For a locally compact space $Z$, let $\mathcal M^1(Z)$ denote its set of Borel probability measures. A family $\mathcal{M}\subset \mathcal{M}^1(Z)$ of probability measures is \emph{uniformly tight} if for every $\varepsilon>0$, some compact set $K\subset Z$ has $\nu(K) \geq 1-\varepsilon$ for every $\nu\in \mathcal{M}$.

\begin{theo}[Uniform tightness of Hecke neighbours]
\label{theo:tightness}
Let $\mu$ be an $A$-invariant, $A$-ergodic probability measure on $X$ with $h_\mu(a_1)>0$. Then $\mathcal T(\mu)$ is uniformly tight. That is, for every $\varepsilon>0$ there is a compact set $K\subset X$ such that $((q_r)_*\lambda)(K)\geq1-\varepsilon$ for every $r\in\comm_G(\Gamma)$ and every good-$r$-lift $\lambda$ of $\mu$.
\end{theo}

\noindent Write $\langle\mathcal T(\mu)\rangle$ for the convex hull of $\mathcal T(\mu)$ in $\mathcal{M}^1(X)$.

\begin{theo}[Positive entropy for ergodic components]
\label{theo:entropy_limits}
Let $\mu$ be as in Theorem~\ref{theo:tightness}. Then every measure $\nu$ in the weak closure of $\langle\mathcal T(\mu)\rangle$ in $\mathcal M^1(X)$ has positive entropy on almost every $A$-ergodic component. More precisely, if $\nu=\int\nu_\xi\,d\eta(\xi)$ is its $A$-ergodic decomposition, then
\[
 h_{\nu_\xi}(a_1)>0\qquad\text{for $\eta$-almost every }\xi.
\]
\end{theo}

There has been a lot of work on non-escape of mass on $X$. If measures are assumed to be invariant under a unipotent element, then by the work of Dani-Margulis \cite{DM93} any weak-$*$ limit of ergodic measures must be the zero measure or a probability measure. This builds on quantitative non-divergence estimates of Dani \cite{Dani1986}, see also \cite{KleinbockMargulis1998}. Stronger rigidity results in this setting were later obtained by Mozes-Shah \cite{MozesShah1995}. 

For general diagonal-invariant measures on $\SL_d(\bbR)/\SL_d(\bbZ)$ there exist weak-$*$ limits $\nu$ where $\nu(X)$ can be anything in the interval $[0,1]$, for example see \cite{Kadyrov2012}. Even for $A$-invariant measures supported on compact $A$-orbits weak-$*$ limits can be wild, see \cite{Shapira2017, SolanYifrach2023}. Under a positive entropy assumption the situation becomes somewhat more rigid. For example the following was shown by Einsiedler, Lindenstrauss, Michel, and Venkatesh.

\begin{theo*}[{\cite[Theorem 5.1]{ELMV2012}}]
    Let $c>0$ and let $(\mu_i)$ be a sequence of $a_1$-invariant probability measures on $X$ with entropy $h_{\mu_i}(a_1)\geq  c$. Then any weak-$*$ limit $\mu$ has mass $\mu(X) \geq c-1$.
\end{theo*}

\subsection{An arithmetic equidistribution theorem}
\label{subsec:intro_padic}

Johnson and Rudolph \cite{JohnsonRudolph1995} proved that if a probability measure on $\bbR/\bbZ$ is invariant and ergodic under multiplication by $r$ and has positive entropy, then its pushforwards under $x\mapsto s^k x$ converge in Ces\`aro average to the Lebesgue measure, provided that $r,s\geq2$ are multiplicatively independent. Einsiedler \cite[Conjecture~9.2]{EinsiedlerICM2010} proposed an analogue for the two diagonal flows on an irreducible quotient of $\SL_2(\bbR)\times\SL_2(\bbR)$, and also asked for versions on other homogeneous spaces that admit some measure rigidity. Our first application verifies an analogue of this conjecture in the $p$-adic extension of $X$, which we now define.

Fix a prime $p$. Set
\begin{equation}
\label{eq:arithmetic_spaces}
 \begin{aligned}
 G_\infty&=\PGL_2(\bbR),&\Gamma_\infty&=\PGL_2(\bbZ),\\
 G_p&=\PGL_2(\bbQ_p),&\Gamma_p&=\PGL_2(\bbZ[1/p]),\\
 K_p&=\PGL_2(\bbZ_p),&Y_p&=(G_\infty\times G_p)/\Gamma_p,
 \end{aligned}
\end{equation}
with $\Gamma_p\leq  G_\infty\times G_p$ embedded diagonally. Dividing by $K_p$ from the left gives the proper $G_\infty$-equivariant map
\begin{equation*}
 \pi:Y_p\longrightarrow K_p\backslash Y_p\simeq X;
 \qquad \pi((g,k)\Gamma_p)=g\Gamma_\infty \text{ for } k\in K_p
\end{equation*}
(see \cite[Section~4.1]{EFS2011}). We let $A$ act on $Y_p$ by multiplication in the real factor, and put
\begin{equation}
\label{eq:padic_element}
b=\left(e,\diag(
    p,1)\right) \in G_\infty \times G_p.
\end{equation}

\begin{theo}
\label{theo:padic_equidistribution}
\label{cor:haar_cesaro}
Let $\mu$ be an $A$-invariant, $A$-ergodic probability measure on $X$ with $h_\mu(a_1)>0$, and let $\sigma\in\mathcal M^1(Y_p)$ be $A$-invariant with $\pi_*\sigma=\mu$. Then
\begin{equation}
\label{eq:haar_cesaro_intro}
 \frac1N\sum_{k=0}^{N-1}\pi_* ((b^k)_*\sigma)\longrightarrow m_X.
\end{equation}
Moreover, there is a set $D\subseteq\bbZ_{\geq0}$ of natural density one such that
\begin{equation}
\label{eq:density_one_measures}
 \pi_*((b^k)_*\sigma)\longrightarrow m_X
 \qquad(k\to\infty,\ k\in D).
\end{equation}
\end{theo}
Here natural density one means that $|D\cap\{0,\ldots,N-1\}|/N\to1$ as $N\to \infty$. 

In a follow-up paper we will use the same methods to prove \cite[Conjecture~9.2]{EinsiedlerICM2010} for an irreducible quotient of $\SL_2(\bbR)\times\SL_2(\bbR)$.

\subsection{Continued fraction patterns}
\label{subsec:intro_patterns}

For an irrational $\alpha$, define the point $x_\alpha:= u^+(-\alpha)\Gamma$ and its forward omega-limit set $\Omega_\alpha$ by
\begin{equation}
\label{eq:omega_definition}
 \Omega_\alpha = \bigcap_{T>0}\overline{\{a_t x_\alpha:t\geq T\}}^{\,X}
 = \{y\in X:a_{t_n}x_\alpha\to y \text{ for some } t_n\to\infty\}.
\end{equation}
The set $\Omega_\alpha$ is closed and $A$-invariant, and is non-empty and compact when $\alpha$ is badly approximable.

For badly approximable $\alpha$, let $P_\alpha(n)$ be the number of distinct patterns of length $n$ occurring in $(a_j(\alpha))_{j\geq1}$. Since the partial quotients are bounded, the limit
\begin{equation*}
 h_{\rm cf}(\alpha):=\lim_{n\to\infty}\frac{\log P_\alpha(n)}{n}
\end{equation*}
exists by subadditivity, and equals the topological entropy of the orbit closure of $\alpha$ under the shift map in its symbolic representation. By the variational principle for this symbolic coding it follows that $h_{\rm cf}(\alpha)>0$ if and only if $\Omega_\alpha$ supports an $A$-invariant probability measure with $h_\mu(a_1)>0$ (see \cite{Walters1982} for the definitions of entropy, and its key properties). The same entropy condition occurs in \cite[Theorem~2.1]{BBEK2015}; see also \cite[Appendix~A]{EFS2011} for the coding.

\begin{theo}
\label{theo:cf_patterns}
Let $\alpha$ be irrational and suppose that $\Omega_\alpha$ supports an $A$-invariant probability measure $\mu$ with $h_\mu(a_1)>0$. For every prime $p$ there is a set $D\subseteq\bbZ_{\geq0}$ of natural density one such that, for every pattern $w$, the pattern $w$ occurs infinitely often in the continued fraction expansion of $p^k\alpha$ for all sufficiently large $k\in D$.

In particular, the conclusion holds if $\alpha$ is badly approximable and $h_{\rm cf}(\alpha)>0$.
\end{theo}

In the badly approximable case, this answers \cite[Question~5.7]{EFS2011} under the stated positive entropy assumption. Note that the set of badly approximable numbers $\alpha$ with $h_{\rm cf}(\alpha)=0$ has Hausdorff dimension zero; see \cite[p.~7]{Bugeaud2014}. Thus the question has a positive answer for every badly approximable number outside a set of Hausdorff dimension zero. Notice that $p^k\alpha$ is still badly approximable, so the minimum threshold for the size of $k$ must depend on $w$.

\begin{rem}
\label{rem:hecke_branches_adelic}
In fact, since Theorems~\ref{theo:tightness} and~\ref{theo:entropy_limits} concern \emph{arbitrary} Hecke neighbours, Theorem~\ref{theo:cf_patterns} holds verbatim for $\{p^{-n}\alpha\}$ instead of $\{p^n\alpha\}$, or more generally for \emph{any} branch in the $p$-Hecke tree based at $x_\alpha$, in the terminology of \cite[\S~4]{AkaShapira2018}. In particular, degenerate branches (see \cite[Def.~4.7]{AkaShapira2018}) play no role here. Indeed, essentially the only change in the argument is in choosing the $p$-adic part of the lift $y_\alpha$ of $x_\alpha$ in~\eqref{eq:lift_choice}, in a way that records the choice of the branch (cf.~\cite[\S~4]{AkaShapira2018}). Furthermore, one can readily extend the argument to $S$-adic spaces and consider arbitrary branches in the $S$-Hecke graph (see again \cite[\S~4]{AkaShapira2018} for the relevant definitions).
\end{rem}

In addition, in selected windows one can obtain the correct statistics. Recall the Gauss map $\Phi:[0,1]\to [0,1]:x \mapsto \frac{1}{x} - \lfloor\frac{1}{x}\rfloor$ (declaring $\Phi(0)=0$) and the Gauss-Kuzmin measure $\tau$ given by \mbox{$d\tau(x) = dx/((1+x)\log 2)$}. Then $\tau$ is a $\Phi$-invariant and ergodic probability measure; see for example \cite[Chapter~3]{EinsiedlerWard2011}. The measure $\tau$ gives statistics of appearances $\tau(I_w)$ of patterns $w$ in the continued fraction expansion of Lebesgue-generic numbers, where $I_w$ is an explicit interval defined by $w$ (cf.~\cite{AkaShapira2018}).   

\begin{cor}[Gauss-Kuzmin statistics in selected windows]
\label{cor:gauss-kuzmin}
    Let $\alpha$ and $D$ be as in Theorem~\ref{theo:cf_patterns}. As $k \to \infty$ in $D$, one can choose positive integers $m_k,\ell_k \to \infty$ such that for any fixed pattern $w$, its frequency in the finite word
    \begin{equation*}
        a_{m_k}(p^k\alpha),\ldots,a_{m_k+\ell_k-1}(p^k\alpha)
    \end{equation*}
    converges as $k\to\infty$ to $\tau(I_w)$.
\end{cor}

We note that quadratic irrationals $\alpha$ have $h_{\rm cf}(\alpha)=0$, and they correspond to closed geodesics in $X$. The uniform measure on a closed geodesic has zero entropy, and using these measures one can give examples showing the necessity of the positive entropy assumption in Theorem~\ref{theo:tightness}.

\begin{exmp}[the golden-ratio geodesic]
    Let $M=\left(\begin{smallmatrix}2&1\\1&1\end{smallmatrix}\right)$ and choose $g\in G$ with $gMg^{-1}=a_T$, where $T=2\log\frac{1+\sqrt5}{2}$. Let $\mu$ be the invariant probability measure on the closed orbit $Ag\Gamma$. Writing $F_j$ for the Fibonacci numbers, set
\[
 r_n=F_{2n}^{-1/2}(e_1,M^ne_1)
        =F_{2n}^{-1/2}\begin{pmatrix}1&F_{2n+1}\\0&F_{2n}\end{pmatrix}.
\]
Notice that $r_n\in\comm_G(\Gamma)$. In fact, since $gM^n = a_{nT}g$, we also see that \mbox{$gr_n = F_{2n}^{-1/2}(v,a_{nT}v)$}, where $v = ge_1$. By the Cayley--Hamilton Theorem we know that \mbox{$M^{2n} = -I+ (\operatorname{Tr} M^n)M^n$}, and so
\begin{equation*}
    M^nr_n = F_{2n}^{-1/2}(M^ne_1, M^{2n}e_1) = r_n \begin{pmatrix}
        0 & -1 \\ 1 & \operatorname{Tr} M^n
    \end{pmatrix}.
\end{equation*}
In particular, $r_n^{-1}M^n r_n \in \Gamma$, and so $M^n \in \Gamma_{r_n}$. Thus $a_{nT}g\Gamma_{r_n} = gM^n \Gamma_{r_n} = g\Gamma_{r_n}$, and so $nT$ is a period of $g\Gamma_{r_n}$, and its invariant probability measure $\lambda_n$ projects to $\mu$. 

The lattice $a_tgr_n\bbZ^2$ contains the non-zero vector $F_{2n}^{-1/2}a_tge_1$, of length $\ll e^{|t|-nT/2}$, where we used that $F_{2n}\asymp e^{nT}$. Therefore for any fixed $\varepsilon$, we have that
\begin{equation*}
   |\{t \in [-\tfrac{nT}{2},\tfrac{nT}{2}]:a_tgr_n\bbZ^2\in X(\varepsilon)\}| = O_\varepsilon(1).
\end{equation*}
Hence $(q_{r_n})_*\lambda_n(X(\varepsilon))=O_\varepsilon(n^{-1})\to0$.
\end{exmp}

Note however that for quadratic irrationals stronger results about the partial quotients of $\{p^k\alpha\}_{k \geq 0}$ and related arithmetic sequences were obtained by the first named author and Uri Shapira in \cite{AkaShapira2018,Aka2020}, using completely different methods.

\subsection{Applications to the \texorpdfstring{$p$}{}-adic Littlewood conjecture}
Applying Theorem~\ref{theo:cf_patterns} to patterns consisting of a single digit gives the following corollary. In the badly approximable case, this recovers \cite[Theorem~2.1]{BBEK2015}.

\begin{cor}
\label{cor:littlewood}
Under the assumptions of Theorem~\ref{theo:cf_patterns}, for every prime $p$,
\[
 \inf_{q\geq1}q\,|q|_p\,\|q\alpha\|=0,
\]
where $|\cdot|_p$ is the $p$-adic absolute value normalised by $|p|_p=p^{-1}$, and $\|\cdot\|$ denotes the distance to the nearest integer.
\end{cor}

Furthermore this answers \cite[Problem~7]{VitorinoVukusic2026}, concerning the set of counterexamples to their $B$-variant of the $2$-adic Littlewood conjecture. For an irrational $\alpha$, write $B(\alpha):=\limsup_{n \to \infty}a_n(\alpha)$. Their question is the case $p=2$ of the following statement.

\begin{cor}
\label{cor:hausdorff_exceptions}
For every prime $p$, the set
\[
 E_p=\left\{\alpha\in\bbR\setminus\bbQ:
             \sup_{k\geq0}B(p^k\alpha)<\infty\right\}
\]
has Hausdorff dimension zero.
\end{cor}
\begin{proof}
If $\alpha\in E_p$, then $B(\alpha)<\infty$, so $\alpha$ is badly approximable. If $h_{\rm cf}(\alpha)>0$, for any positive integer $m$ we apply Theorem~\ref{theo:cf_patterns} to the one-digit pattern $m$; this gives an exponent $k$ with $B(p^k\alpha)\geq m$. This contradicts $\alpha\in E_p$. Hence
\[
 E_p\subseteq\{\alpha\text{ badly approximable}:h_{\rm cf}(\alpha)=0\}.
\]
It is a standard exercise to show that the set on the right-hand side has Hausdorff dimension zero; see \cite[Proposition~1]{MauduitMoreira2012}.
\end{proof}

\subsection{Overview of the argument and organization of the paper}
\label{subsec:intro_overview}

Theorem~\ref{theo:padic_equidistribution} follows from Theorems~\ref{theo:tightness} and~\ref{theo:entropy_limits} using \cite[Theorem~1.1]{Lindenstrauss2006}; see Section~\ref{subsec:haar_limits}. Theorem~\ref{theo:cf_patterns} then follows from Theorem~\ref{theo:padic_equidistribution} by the standard connection between the geodesic flow and continued fraction expansions; see Section~\ref{subsec:cf_coding}. We therefore give an overview of the proof of Theorem~\ref{theo:tightness}.

To prove Theorem~\ref{theo:tightness}, we first want to couple the ``source of randomness'', i.e.\ the positive-entropy measure $\mu$, with its Hecke neighbours. We therefore define, for a good-$r$-lift $\lambda\in\mathcal M^1(X_r)$ with $(p_r)_*\lambda=\mu$, the measure
\[
 J_{r,\lambda}=(p_r,q_r)_*\lambda\in\mathcal M^1(X\times X).
\]
Note that its first marginal is $\mu$ and its second marginal is a Hecke neighbour of $\mu$, for which we want to prove uniform tightness.

A feature of our approach is to study escape of mass by taking the closure $\overline X$ of $X$ in $\operatorname{Sub}(\bbR^2)$, the space of closed subgroups of $\bbR^2$ with the Chabauty topology. Its boundary consists precisely of the one-dimensional vector subspaces of $\bbR^2$, and $\overline X$ admits the following decomposition
\[
 \overline X=X\sqcup\bbP^1(\bbR)
\]
as a $G$-space (see Section~\ref{subsec:compactification_by_lines}). A weak-$*$ limit $J$ of the $J_{r,\lambda}$ is a probability measure on $X\times\overline X$, since the first projection remains $\mu$. To prove tightness, it suffices to show that $J(X\times\bbP^1(\bbR))=0$. By the $A$-invariance of $J$, any escaped mass is supported on the two fixed directions $[e_1]$ and $[e_2]$ in $\bbP^1(\bbR)$. Here $e_1$ and $e_2$ are the standard basis vectors of $\bbR^2$.

To show that mass also cannot escape to $[e_1]$ or $[e_2]$, we define in Sections~\ref{subsec:probability_kernels} and \ref{subsec:leafwise_measures} averaging operators on $\mathcal M^1(X)$ using the $U^\pm$-leafwise measures of $\mu$, these averaging operators are non-trivial by the positive-entropy assumption on $\mu$. In Section~\ref{subsec:hecke_neighbours} we extend these averaging operators from $\mathcal M^1(X)$ to $\mathcal M^1(X_r)$ and then to $\mathcal M^1(X\times X)$ in a compatible way that keeps track of the source of randomness $\mu$. We then prove that such averaging commutes with weak-$*$ limits, mainly because the averaging is with respect to the fixed first marginal (see Lemma~\ref{lem:limits_of_joinings}). 

Since positive entropy implies that the averaging procedure is non-trivial, and every element of $U^+\setminus\{e\}$ moves $[e_2]$ away from $\{[e_1],[e_2]\}$, it will follow that $[e_2]$ carries no mass (see the calculation \eqref{eq:no_escape_of_mass}). A similar argument will apply with $U^-$ and $[e_1]$. The desired uniform tightness follows.

\subsection*{Acknowledgments}
The first author thanks Uri Shapira for introducing him to homogeneous dynamics through the problem {\cite[Question~5.7]{EFS2011}} stated above. The authors thank Manfred Einsiedler and Elon Lindenstrauss for their interest and for improving earlier versions of our arguments.

\subsection*{AI declaration}

The authors acknowledge the use of OpenAI's GPT-6-Astra in developing some of the arguments, accessed through a commercial subscription. The key elements in the proofs are the authors' own. The text in this article is entirely human-written, and the authors take sole responsibility for any mistakes in this work.

    \section{Tightness and entropy of Hecke neighbours}
\label{sec:hecke_dynamics}

\noindent In this section we prove Theorems~\ref{theo:tightness} and~\ref{theo:entropy_limits}. We use the notation of Section~\ref{subsec:intro_hecke}; throughout $\mu$ is an $A$-invariant probability measure on $X$. We will later assume in addition that it is $A$-ergodic with $h_\mu(a_1)>0$.

\subsection{Probability kernels and averaging}
\label{subsec:probability_kernels}

Let $I = [-1,1]$ and let $Z$ be a second countable locally compact space.

\begin{defn}
    A \emph{(Borel) probability kernel} on $Z$ is a measurable function $\kappa:Z \to\mathcal{M}^1(I)$. We denote $\kappa(z)$ by $\kappa_z$.
\end{defn}

The probability kernel is measurable in the sense that for every measurable $B \subset I$, the map $z \mapsto \kappa_z(B)$ is measurable. Suppose now that $Z$ is equipped with a continuous action of a group $U \cong \bbR$, parameterised as $U = \{u(s):s\in \bbR\}$, and let $\nu$ be any probability measure on $Z$. Then we define the averaged measure $\kappa *\nu$ by setting for every $f \in C_c(Z)$
\begin{equation*}
    (\kappa * \nu)(f):= \int_Z \int_I f(u(s)z)\, d\kappa_z(s)\, d\nu(z).
\end{equation*}

\subsection{Leafwise measures}
\label{subsec:leafwise_measures}

Notice that by Poincar\'{e} recurrence, the $U^\pm$-action on $X$ is free $\mu$-almost everywhere. The following is a well-known construction, see \cite[Section~3]{Fubini} for a modern treatment. Recall that for any two measures $\tau_1,\tau_2$ we write $\tau_1 \propto \tau_2$ if there is some $c>0$ such that $\tau_1 = c\tau_2$.

\begin{prop}[{\cite[Theorem 6.3]{Pisa}}]
\label{prop:leafwise_measures}
Let $U = U^+$ or $U^-$. There is a conull set $X_0 \subset X$ such that for every $x \in X_0$ there is a Radon measure $\mu_x^U$ on $\bbR$ called a \emph{leafwise measure}, satisfying:
    \begin{enumerate}[label = {(\roman*)}]
        \item For every $f\in C_c(\bbR)$, $x \mapsto \int_\bbR f \, d\mu_x^U$ is measurable; 
        \vspace{1mm}
        \item For every $x\in X_0$ we have that $0 \in\supp (\mu_x^U)$;
        \vspace{1mm}
        \item \label{prop:leafwise_measures:compatibility} For every $x\in X_0$ and $u\in U$ such that $ux\in X_0$ we have that $\mu_x^U \propto u_*(\mu_{ux}^U)$;
        \vspace{1mm}
        \item Let $B\subset X$ and let $\mathcal{A}$ be any countably generated $\sigma$-subalgebra of subsets of $B$ such that for any $x\in B$, the atom $[x]_\mathcal{A}$ is an \emph{open $U$-plaque} (that is, $U_{x,\mathcal{A}}:=\{u\in U:ux\in [x]_\mathcal{A}\}$ is open and bounded in $U$). Then for $\mu$-almost every $x\in B$,
        \begin{equation*}
            (\mu|_B)_x^\mathcal{A} \propto (\mu_x^U|_{U_{x,\mathcal{A}}})\cdot x
        \end{equation*}
        (the right hand side denotes the pushforward under the orbit map at $x$).
    \end{enumerate}
    These properties define $\mu_x^U$ on $X_0$ uniquely up to scaling.
\end{prop}

\noindent We will use $\mu^\pm_x$ to denote the measures $\mu^{U^\pm}_x$. We choose to normalise these leafwise measures such that $\mu_x^\pm(I) = 1$. We can now define probability kernels on $X$ by restricting the leafwise measure to $I$. That is, we set $\kappa_x^\pm := \mu_x^\pm|_{I}$ for $x\in X_0$ and $\kappa_x^\pm = \delta_0$ otherwise. 

Notice that for each $\mu$-null set $N \subset X$, $\mu_x^\pm(\{s:u^\pm(s)x\in N\}) = 0$ for $\mu$-almost every $x$. It follows that $\kappa^\pm *\mu(N) = 0$, and so in particular $\kappa^\pm *\mu \ll \mu$. Define
\begin{equation*}
    \phi^\pm(x) := \frac{d(\kappa^\pm*\mu)}{d\mu}(x)
\end{equation*}
to be the Radon-Nikodym derivative, so that $\int \phi^\pm \, d\mu = 1$, and as measures
\begin{equation}
\label{eq:averaging_on_X}
    d(\kappa^\pm *\mu) = \phi^\pm d\mu.
\end{equation} 
The $A$-ergodicity and positive entropy of $\mu$ tell us that these probability kernels are almost everywhere non-trivial.

\begin{prop}[{\cite[Theorem~7.6~(iii)]{Pisa}}]
\label{prop:positive_entropy}
    If $\mu$ is in addition $A$-ergodic and \mbox{$h_\mu(a_1)>0$}, then $\kappa_x^\pm(\{0\}) = 0$ for $\mu$-almost every $x$. 
\end{prop}

It is easy to see that under the assumptions of Proposition~\ref{prop:positive_entropy}, in fact $\kappa_x^\pm$ is non-atomic almost everywhere.

\subsection{Hecke neighbours}
\label{subsec:hecke_neighbours}

We want to understand to what extent \eqref{eq:averaging_on_X} extends to finite covers of $X$ and to Hecke neighbours. We begin with the former.

\begin{lem}
\label{lem:averaging_the_lift}
    Let $\Delta \leq \Gamma$ be finite index, and consider $Y = G/\Delta$. Let $p:Y \to X$ be the natural projection. Let $\lambda \in \mathcal{M}^1(Y)$ be an $A$-invariant probability measure such that $p_*\lambda = \mu$. Define the probability kernel $\hat{\kappa}^\pm:Y \to \mathcal{M}^1(I)$ by $\hat{\kappa}^\pm_y := \kappa^\pm_{p(y)}$. Then as measures
    \begin{equation*}
        d(\hat{\kappa}^\pm*\lambda) = (\phi^\pm \circ p)\,  d\lambda. 
    \end{equation*}
\end{lem}

\begin{proof}
    The proof works in two steps. Firstly we define a uniform lift of $\mu$ and show the result there, and then compare the uniform lift to $\lambda$.

    For any Borel function $f$ on $Y$, define a function $f^p$ on $X$ via $f^p(x) := \frac{1}{m}\sum_{p(y)=x}f(y)$, where $m = [\Gamma:\Delta]$. Let $\hat{\mu}$ be the measure on $Y$ given by $\hat{\mu}(f) := \mu(f^p)$, this is the \emph{uniform lift of $\mu$}. This measure is clearly $A$-invariant, and we can compute for any $f\in C_c(Y)$
    \begin{align*}
        (\hat{\kappa}^\pm*\hat{\mu})(f) & = \int_Y\int_I f(u(s)y)\, d\hat{\kappa}^\pm_y(s)\, d\hat{\mu}(y) \\
        & = \int_X\int_I\frac{1}{m}\sum_{p(y)=x}f(u(s)y)\, d\kappa^\pm_x(s)\, d\mu(x) \\
        & = \int_X\int_I f^p(u(s)x)\, d\kappa_x^\pm(s)\, d\mu(x)\\
        & = \int f^p(x) \phi^\pm(x)\, d\mu(x) \\
        & = \int_Y f(y)(\phi^\pm\circ p)(y)\, d\hat{\mu}(y) 
    \end{align*}
    where we used that $p$ is $U^\pm$-equivariant in the third line, and that $(\phi^\pm \circ p)$ is constant on each fiber of $p$ in the last line. So in particular $d(\hat{\kappa}^\pm*\hat{\mu}) = ( \phi^\pm\circ p)\, d\hat{\mu}$. 

    Consider now any Borel set $E \subset Y$. Then $\lambda(E) \leq \mu(p(E)) \leq m\cdot  \hat{\mu}(E)$ and hence $\lambda \ll m \hat{\mu}$ and so $d\lambda = \varphi\, d\hat{\mu}$ where $\varphi\in L^1(\hat{\mu})$ is the Radon-Nikodym derivative. Since both $\lambda$ and $\hat{\mu}$ are $A$-invariant, we have that $\varphi \circ a_1 = \varphi$ $\hat{\mu}$-almost everywhere. Therefore by the Hopf argument (for example, \cite[Lemma~7.20]{Pisa} applied countably many times) we can find a function $\varphi^+$ that agrees with $\varphi$ $\hat{\mu}$-almost everywhere, and that is in addition $U^+$-invariant. We thus calculate
    \begin{align*}
        (\hat{\kappa}^+*\lambda)(f)& = \int_Y \int_I f(u(s)y)\, d\hat{\kappa}^+_y(s)\, d\lambda(y) \\
        & = \int_Y \int_I f(u(s)y)\, d\hat{\kappa}^+_y(s)\, \varphi^+(y) d\hat{\mu}(y) \\
        & = \int_Y \int_I f(u(s)y) \varphi^+(u(s)y)\, d\hat{\kappa}^+_y(s)\,  d\hat{\mu}(y) \\
        & = \int_Y f(y) \varphi^+(y) (\phi^+ \circ p)(y)\, d\hat{\mu}(y) \\
        & =\int_Y f(y)(\phi^+ \circ p)(y)\, d\lambda(y).
    \end{align*}
    Repeating this argument with $\hat{\kappa}^-$ and $U^-$ gives the claim.
\end{proof}

\begin{rem}
    It was pointed out to the authors by Elon Lindenstrauss that this argument could be made more conceptual by directly applying \cite[Proposition 2.10]{EinsiedlerLindenstrauss2018} to $\lambda$, without considering $\hat{\mu}$ first. Indeed in our setting \cite[Proposition 2.10]{EinsiedlerLindenstrauss2018} says that $\lambda_y^U = \mu_{p(y)}^U$ for $\lambda$-almost every $y$, which is also shown using an argument reminiscent of the Hopf argument. It is also somewhat analogous to the product structure argument developed in \cite{EinsiedlerKatok2005}.
\end{rem}

We are now ready to transfer the identity \eqref{eq:averaging_on_X} to the measures $J_{r,\lambda}$. Recall that $J_{r,\lambda}:= (p_r,q_r)_*\lambda \in \mathcal{M}^1(X\times X)$, for some $r\in \comm_G(\Gamma)$ and $\lambda$ a good-$r$-lift.

\begin{prop}
\label{prop:averaging_the_joining}
    For every $r\in \comm_G(\Gamma)$ and good-$r$-lift $\lambda$, we have that
    \begin{equation*}
        d(\kappa_1^\pm * J_{r,\lambda}) = (\phi^\pm\circ\pi_1) \, d J_{r,\lambda},
    \end{equation*}
    where the averaging is taken under the diagonal $(u(s),u(s))$ action on $X\times X$ for the probability kernels $\kappa^\pm_1 := \kappa^\pm \circ \pi_1$. Here $\pi_1:X\times X \to X$ denotes the projection onto the first coordinate.
\end{prop}

\begin{proof}
        By Lemma~\ref{lem:averaging_the_lift}, for any good-$r$-lift $\lambda$ on $X_r$ we have $d(\hat{\kappa}^\pm *\lambda) = (\phi^\pm \circ p_r)\,  d\lambda$. Since $(p_r,q_r)$ is $G$-equivariant and $\hat{\kappa}_y^\pm = \kappa_{p_r(y)}^\pm$, we push this forward to obtain the result.
\end{proof}


\subsection{Compactification by lines}
\label{subsec:compactification_by_lines}


Let $\pi_i:X\times X \to X$ be the projection to the $i^{th}$ coordinate. We want to be able to take limits of the $J_{r,\lambda}$'s, so we need to compactify $X \times X$. Since by definition $(\pi_1)_*J_{r,\lambda} = \mu$ always, it suffices to compactify the second factor. 

To do this, consider the set $\operatorname{Sub}(\bbR^2)$ of closed subgroups of $\bbR^2$ equipped with the Chabauty topology. In particular, a sequence $H_n \in \operatorname{Sub}(\bbR^2)$ converges to some $H \in \operatorname{Sub}(\bbR^2)$ if and only if:
\begin{enumerate}
    \item For every subsequence $h_k \in H_{n_k}$ such that $h_k \to h \in \bbR^2$, in fact $h\in H$, and
    \item Every $h \in H$ is the limit in $\bbR^2$ of some $h_n \in H_n$.
\end{enumerate}
The space $\operatorname{Sub}(\bbR^2)$ with this topology is compact, and admits a continuous left $G$-action. We can view $X \subset \operatorname{Sub}(\bbR^2)$ and let $\overline{X}$ be the closure of $X$ in $\operatorname{Sub}(\bbR^2)$, which is compact.

 Consider a sequence of lattices $\Lambda_n \in X$ that converge to some $H \in \operatorname{Sub}(\bbR^2)$. Either they converge in $X$, or the shortest non-zero vector $v_n \in \Lambda_n \setminus \{0\}$ has $\varepsilon_n:= \|v_n\| \to 0$. In this case, in particular $v_n$ is primitive, and so we can extend it to an integral basis $(v_n,w_n)$ for $\Lambda_n$. Now notice that $d(w_n,\bbR v_n) = \frac{|\det(v_n,w_n)|}{\|v_n\|} = \frac{1}{\varepsilon_n}$, and by the definition of the Chabauty topology we must have that $\Lambda_n \to \bbR v$, where $\bbR v_n \to \bbR v$. Conversely, every subspace $\bbR v \in \operatorname{Sub}(\bbR^2)$ is a limit of a sequence of unimodular lattices. Thus the closure of $X$ in $\operatorname{Sub}(\bbR^2)$ can be identified as $\overline{X} = X\sqcup \bbP^1(\bbR)$. We collect these useful facts in the following statement, for a full treatment on the topology of $\overline{X}$ see \cite[Theorem~4.6~(i)]{BHK2009}.

\begin{lem}
    The space $X$ is open in $\overline{X}$, and the subspace topology coincides with the topology coming from $G/\Gamma$. In addition, $\overline{X} \setminus X \cong \bbP^1(\bbR)$, both $X$ and $\bbP^1(\bbR)$ are $G$-invariant, and the action of $G$ on $\bbP^1(\bbR)$ is the usual projective action.
\end{lem}

The point in $\bbP^1(\bbR)$ corresponding to the subspace $\bbR v$ will be denoted by $[v]$. We view any measure on $X\times X$ as one on $X\times \overline{X}$, which will allow us to take weak-$*$ limits. The relation in Proposition~\ref{prop:averaging_the_joining} allows us to control these.

\begin{lem}
\label{lem:limits_of_joinings}
    Let $J_n \in \mathcal{M}^1(X\times \overline{X})$ be a sequence such that $(\pi_1)_*J_n = \mu$ for all $n$. Suppose that $J_n \to J$ in the weak-$*$ topology. Then $J$ is a probability measure on $X\times \overline{X}$. If in addition we have a probability kernel $\kappa$ on $X$ and a function $\phi \in L^1(\mu)$, then as measures
    \begin{equation*}
        \kappa_1*J_n \to \kappa_1*J \quad \text{and} \quad (\phi \circ \pi_1) \cdot J_n \to (\phi\circ\pi_1) \cdot J,
    \end{equation*}
    where $\kappa_1 = \kappa \circ\pi_1$. In particular, if $d(\kappa_1*J_n) = (\phi\circ\pi_1)\, d J_n$ for all $n$, then $d(\kappa_1*J) = (\phi\circ\pi_1)\, d J$ as measures.
\end{lem}

\begin{proof}
    It is clear that $(\pi_1)_*J=\mu$ and that $J$ is a probability measure. By Lusin's Theorem for any $\varepsilon>0$ there is a compact set $K\subset X$ with $\mu(X\setminus K) < \varepsilon$ on which $x\mapsto \kappa_x$ is continuous (for the weak topology on $\mathcal{M}^1(I)$). Consider any $f \in C_c(X\times \overline{X})$, and consider 
    \begin{equation*}
        \kappa_1*f(x,y) := \int f(u(s)x,u(s)y)\, d\kappa_x(s),
    \end{equation*}
    we have that $(\kappa_1*f)|_{K\times \overline{X}}$ is continuous. By Tietze's Extension Theorem we can therefore extend $\kappa_1 * f$ to a function $h \in C_c(X\times \overline{X})$ that satisfies $\|h\|_\infty \leq \|f\|_\infty$. For $J' = J_n$ or $J$ we have 
    \begin{equation*}
        |J'(\kappa_1*f - h)| \leq 2 \|f\|_\infty \mu(X\setminus K)\leq 2\varepsilon \|f\|_\infty,
    \end{equation*}
    where we used that the first marginal is always $\mu$. Since $J_n(h) \to J(h)$ by weak convergence, letting $\varepsilon \to 0$ we also see that $\kappa_1 *J_n(f) \to \kappa_1*J(f)$.
    For the second convergence, choose $s_k\in C_c(X)$ such that $\|s_k - \phi\|_{L^1(\mu)} \to 0$. Then uniformly for $J' = J_n$ or $J$ we have that
    \begin{equation*}
        |J'(((\phi-s_k) \circ \pi_1)f)| \leq \|f\|_\infty \|\phi-s_k\|_{L^1(\mu)}.
    \end{equation*}
    Since $(s_k\circ\pi_1 )f \in C_c(X\times \overline{X})$ we have that $J_n((s_k\circ\pi_1 )f) \to J((s_k\circ\pi_1 )f)$, and the result follows.
\end{proof}

The following is an immediate corollary of Proposition~\ref{prop:averaging_the_joining} and Lemma~\ref{lem:limits_of_joinings}.

\begin{cor}
\label{cor:limit_is_abs_cont}
    Let $J$ be any weak-$*$ limit of some $J_{r,\lambda}$'s. Then $\kappa^{\pm }_1 * J \ll J$. 
\end{cor}

We are now ready to prove Theorem~\ref{theo:tightness}, by considering what mass weak-$*$ limits could have on the boundary.

\begin{proof}[Proof of Theorem~\ref{theo:tightness}]
    Let $J$ be any weak-$*$ limit of some measures $J_{r,\lambda}$, and notice in particular that it is $A$-invariant. Consider $(\pi_2)_*(J|_{X\times \bbP^1(\bbR)})$, this is an $A$-invariant measure on $\bbP^1(\bbR)$ and so has support in $\{[e_1], [e_2]\}$. Let $E = X\times (\bbP^1(\bbR)\setminus \{[e_1],[e_2]\})$, so we have that $J(E) = 0$. By Corollary~\ref{cor:limit_is_abs_cont} we also have that $\kappa_1^\pm *J(E) = 0$. 

    We consider now the $U^+$-action on $\bbP^1(\bbR)$. Notice that $u^+(s)[e_1] = [e_1]$ for all $s$, and that $u^+(s)[e_2] \in  \{[e_1],[e_2]\}$ if and only if $s=0$. We can therefore calculate
    \begin{align}
    \label{eq:no_escape_of_mass}
        0 = (\kappa_1^+*J)(E) & = \int_{X \times \bbP^1(\bbR)} \int_{I} \chi\{(u(s)x,u(s)[v])\in E\}\, d\kappa_x^+(s)\, dJ(x,[v]) \nonumber\\
        & = \int_{X \times \bbP^1(\bbR)} \int_{I} \chi\{u(s)[v] \notin \{[e_1],[e_2]\}\}\, d\kappa_x^+(s)\, dJ(x,[v]) \nonumber\\
        & = \int_{X \times \{[e_2]\}}(1-\kappa_x^+(\{0\}))\, dJ(x,y) \nonumber\\
        & = J(X \times \{[e_2]\}).
    \end{align}
   Here we used Proposition~\ref{prop:positive_entropy} to note that $\kappa^+_x(\{0\}) = 0$ for $\mu$-almost every $x$, and in the third line we used that $E$ is a $J$-null set. Similarly by considering $U^-$ we see that $J(X \times \bbP^1(\bbR)) =0$, and so $J$ is really a probability measure on $X\times X$, as required. By Prokhorov's Theorem we conclude that $\mathcal{T}(\mu)$ is uniformly tight.

    To spell this out: If $\mathcal{T}(\mu)$ were not uniformly tight, there would be an $\varepsilon>0$ and $\nu_n \in \mathcal{T}(\mu)$ with $\nu_n(X(\frac{1}{n})) < 1-\varepsilon$. Consider the associated joinings and (by taking a subsequence) suppose that $J_n \to J$ as above, and let $\nu = (\pi_2)_*J \in \mathcal{M}^1(X)$. Choose some $f \in C_c(X)$ with $0\leq f\leq 1$ and $\nu(f) > 1-\varepsilon/2$, viewed as a function in $C(\overline{X})$. Since $\supp f$ is compact we know that for large enough $n $, $\supp f \subset X(\frac{1}{n})$ but we assumed that $\nu_n(f) < 1-\varepsilon$, contradicting $\nu_n \to \nu$.
\end{proof}

\subsection{Positive entropy in the limiting projections}
\label{subsec:componentwise_positive_entropy}

The following statement is standard and follows from the compatibility of leafwise measures under translation by $U$ (Proposition~\ref{prop:leafwise_measures}~\ref{prop:leafwise_measures:compatibility}) and from \cite[Theorem~7.6~(ii)--(iii)]{Pisa}.

\begin{lem}
\label{lem:zero_entropy_measures}
Let $\nu$ be an $A$-invariant probability measure on $X$ such that $h_\nu(a_1)=0$. Then there is a subset $X_\nu\subset X$ of full $\nu$-measure such that for every $x \in X$, the set $\{s\in \bbR:u(s)x\in X_\nu\}$ has at most one element.
\end{lem}

We are now ready to prove Theorem~\ref{theo:entropy_limits}.

\begin{proof}[Proof of Theorem~\ref{theo:entropy_limits}]
Let $\nu$ be a measure in the weak-$*$ closure of $\langle \mathcal{T}(\mu)\rangle$, and choose finite convex combinations $\nu_n$ of elements of $\mathcal T(\mu)$ that converge weakly to $\nu$. Taking the same combinations of the associated Hecke joinings gives measures $J_n$ such that $(\pi_1)_*J_n = \mu$ and $(\pi_2)_*J_n = \nu_n$. Up to passing to a subsequence we can assume that $J_n \to J$ weakly, where $(\pi_1)_*J = \mu$ and $(\pi_2)_*J = \nu$. Recall that by Proposition~\ref{prop:averaging_the_joining} and Lemma~\ref{lem:limits_of_joinings} we have that $d(\kappa_1^+*J) = (\phi^+\circ \pi_1)\, dJ$.

Let $Z$ be the union of the zero-entropy components in the $A$-ergodic decomposition of $\nu$ and put $c=\nu(Z)$. Suppose that $c>0$. Set $\nu_0 := c^{-1}\nu|_Z$. Since entropy is affine for the ergodic decomposition this is an $A$-invariant probability measure such that $h_{\nu_0}(a_1) = 0$. By Lemma~\ref{lem:zero_entropy_measures}, there is a $\nu_0$-conull Borel set $Z_0\subset Z$ such that for any $x\in X$ the set $\{s:u^+(s)x\in Z_0\}$ has at most one element. We can now calculate from the definition of the averaging that
\begin{align*}
 (\kappa_1^+*J)(X\times Z_0)
 &=\int_{X\times X}\kappa^+_x(\{s\in I:u^+(s)y\in Z_0\})\,dJ(x,y)\\
 &=0.
\end{align*}
We used here that the $\kappa_x^+$ are non-atomic since $\mu$ has positive entropy (Proposition~\ref{prop:positive_entropy}). Notice that $(\pi_1)_*(J|_{X\times Z})$ is an $A$-invariant measure dominated by $\mu$. Since $\mu$ was assumed to be $A$-ergodic, it follows that $(\pi_1)_*(J|_{X\times Z}) = c\mu$. As such, we can compute
\begin{align*}
    (\kappa_1^+ *J)(X\times Z_0) = \int_{X\times Z_0}(\phi^+ \circ \pi_1)\, dJ = c\int_X\phi^+ \, d\mu = c,
\end{align*}
and so $c = 0$.
\end{proof}

    \section{Arithmetic applications}

We first prove Theorem~\ref{cor:haar_cesaro}, using the two results of Section~\ref{sec:hecke_dynamics}. Recall the notation introduced in \eqref{eq:arithmetic_spaces} and \eqref{eq:padic_element}.

\subsection{The arithmetic extension and Hecke correspondences} 

We recall how to relate pushforwards of measures under $b^k$ to Hecke neighbours. This is immediate once the necessary identifications are made.

\begin{lem}
\label{lem:hecke_realization}
Let $\sigma$ be an $A$-invariant probability measure on $Y_p$ such that $\pi_*\sigma=\mu$. Then for every $k\geq 0$ we have that $\pi_*((b^k)_*\sigma)\in\mathcal T(\mu)$.
\end{lem}

\begin{proof}
Let $d = \diag(p,1)$ and $r_k = p^{k/2}d^{-k}$. We will construct a good-$r_k$-lift $\lambda$ on $X_{r_k}$, such that $(q_{r_k})_*\lambda = \pi_*((b^k)_*\sigma)$. To do this, set $H_k = K_p \cap d^{-k}K_p d^k$, and recall the identification $\pi:Y_p \to K_p\backslash Y_p \cong X$. Clearly both $y\mapsto \pi(y)$ and $y \mapsto \pi(b^k y)$ are constant on $H_k$-orbits.

By strong approximation for $\SL_2$ we have that $G_p = H_k \Gamma_p$, and thus every $H_k$-orbit in $Y_p$ has a representative $(g,e)\Gamma_p$ for some $g\in G_\infty$. Two such representatives $(g,e)\Gamma_p$ and $(g',e)\Gamma_p$ give the same $H_k$-orbit exactly when $g^{-1}g' \in \Gamma_p \cap H_k = \Gamma_\infty \cap d^{-k}\Gamma_\infty d^k =: \Delta_k$. Hence $H_k \backslash Y_p \cong G_\infty / \Delta_k$.

Since $\diag(-1,1) \in \Delta_k$, in fact every coset has a representative in $G$. The inverse image of $\Delta_k\cap\PSL_2(\bbZ)$ in $G$ is $\Gamma_{r_k}$ and hence in fact $H_k \backslash Y_p \cong X_{r_k}$ via the identification $H_k(g,e)\Gamma_p  \leftrightarrow g\Gamma_{r_k}$. Denote the resulting projection map by $\rho_k: Y_p \to X_{r_k}$.

Now notice that $(g,d^k)\Gamma_p = (gd^{-k},e)\Gamma_p = (gr_k,e)\Gamma_p$, and so $q_{r_k} \circ \rho_k = \pi \circ b^k$. Setting $\lambda = (\rho_k)_*\sigma$ the result follows.
\end{proof}

\subsection{Equidistribution}
\label{subsec:haar_limits}

\begin{proof}[Proof of Theorem~\ref{cor:haar_cesaro}]
Let $\sigma_k=(b^k)_*\sigma$ and $\theta_N=N^{-1}\sum_{k=0}^{N-1}\sigma_k$. By Lemma~\ref{lem:hecke_realization}, we have that $\pi_*\sigma_k\in\mathcal T(\mu)$, so
\[
 \pi_*\theta_N=\frac1N\sum_{k=0}^{N-1}\pi_*\sigma_k
 \in\langle\mathcal T(\mu)\rangle.
\]
Theorem~\ref{theo:tightness} and the fact that $\pi$ is proper imply that the $\theta_N$ are uniformly tight. Let $\theta$ be any weak-$*$ limit of $(\theta_N)$, and note that it is $A$-invariant since all the $\sigma_k$'s are. By a telescoping argument it is $b$-invariant: $b_*\theta_N-\theta_N=N^{-1}(\sigma_N-\sigma)$. Its projection $\nu:=\pi_*\theta$ lies in the weak closure of $\langle\mathcal T(\mu)\rangle$ in $\mathcal{M}^1(X)$.

By Theorem~\ref{theo:entropy_limits}, almost every $A$-ergodic component of $\nu$ has positive entropy for $a_1$. By the uniqueness of the ergodic decomposition, almost every $A$-ergodic component of $\theta$ projects to an $A$-ergodic component of $\nu$. Since entropy decreases under factor maps, almost every $A$-ergodic component of $\theta$ also has positive entropy. We are now in a position to apply \cite[Theorem~1.1]{Lindenstrauss2006}, which says in particular that $\theta$ is $\SL_2(\bbR)$-invariant; and hence $\pi_*\theta=m_X$. This proves \eqref{eq:haar_cesaro_intro}.


We now prove \eqref{eq:density_one_measures}, which is a standard consequence of convexity. Set $\nu_k=\pi_*\sigma_k$ and let $\mathcal K$ be the weak closure of $\{\nu_k:k\geq0\}$ in $\mathcal M^1(X)$. By Theorem~\ref{theo:tightness}, $\mathcal K$ is a compact metrisable set of $A$-invariant probability measures. Consider the measure
\[
 \eta_N:= \frac1N\sum_{k=0}^{N-1}\delta_{\nu_k}\in\mathcal M^1(\mathcal K).
\]
If $\eta_N \to \eta$ along a subsequence, then \eqref{eq:haar_cesaro_intro} says that 
\begin{equation*}
    \int_\mathcal{K} \rho \, d\eta(\rho) = m_X.
\end{equation*}
This expresses $m_X$ as a mixture of $A$-invariant probability measures. Since $m_X$ is $A$-ergodic it follows that in fact $\eta = \delta_{m_X}$. In particular, for any fixed neighbourhood $U$ of $m_X$ in $\mathcal{M}^1(X)$, the set $\{k:\nu_k \notin U\}$ has density zero. We can now define $D$ via a diagonal argument. Choose a neighbourhood basis $U_1 \supset U_2 \supset \cdots $ of open sets at $m_X$, and choose integers $N_1 < N_2 < \cdots$ such that 
\begin{equation*}
    \#\{0 \leq k < N: \nu_k \notin U_j\} \leq 2^{-j}N \text{ for } N \geq N_j.
\end{equation*}
Now let $D$ contain $0,\ldots, N_1-1$, and for $N_j \leq k < N_{j+1}$ include $k$ in $D$ precisely when $\nu_k \in U_j$.

$D$ has density one. Indeed if $N_j \leq N < N_{j+1}$, every $ k < N$, $k\notin D$ satisfies $\nu_k \notin U_i$ for some $i \leq j$, and in particular $\nu_k \notin U_j$. By choice of $N_j$ there are at most $2^{-j}N$ of these. Finally, if $k\in D$ and $k> N_j$ then $\nu_k \in U_j$; this gives convergence to $m_X$ along $D$.
\end{proof}


\subsection{From equidistribution to continued fraction patterns}
\label{subsec:cf_coding}
\label{subsec:matrix_calculation}

We now deduce Theorem~\ref{theo:cf_patterns} from Theorem~\ref{theo:padic_equidistribution}. We first construct an invariant lift of the measure supported on $\Omega_\alpha$, then relate translation by $b$ to multiplication of $\alpha$ by $p$, and finally use continued fraction coding to obtain the required recurrence of patterns.

Fix $\alpha$ as in Theorem~\ref{theo:cf_patterns}. Since entropy is affine for the ergodic decomposition, we can assume that $\mu$ is an $A$-invariant, $A$-ergodic probability measure supported on $\Omega_\alpha$, with $h_\mu(a_1)>0$. Consider 
\begin{equation}
\label{eq:lift_choice}
 y_\alpha=(u^+(-\alpha),e)\Gamma_p\in Y_p,\qquad
 \widetilde\Omega_\alpha=\bigcap_{S>0}\overline{\{a_ty_\alpha:t\geq S\}}^{\,Y_p}.
\end{equation}
For every irrational $\alpha$, the set $\widetilde\Omega_\alpha$ is closed and $A$-invariant. Properness of $\pi$ gives
\begin{equation}
\label{eq:omega_projection}
 \pi(\widetilde\Omega_\alpha)=\Omega_\alpha.
\end{equation}
To apply Theorem~\ref{theo:padic_equidistribution}, we need to lift $\mu$ to an $A$-invariant measure on $\widetilde\Omega_\alpha$.

\begin{lem}
\label{lem:supported_lift}
There is an $A$-invariant probability measure $\sigma$ supported on $\widetilde\Omega_\alpha$ such that $\pi_*\sigma=\mu$.
\end{lem}

\begin{proof}
Choose a Borel section $s:\Omega_\alpha \to \widetilde\Omega_\alpha$ of the proper surjection $\pi:\widetilde\Omega_\alpha \to \Omega_\alpha$, and let $\lambda = s_*\mu$. For $T>0$ set
\[
 \lambda_T=\frac1T\int_0^T(a_t)_*\lambda\,dt.
\]
These measures are supported on $\widetilde\Omega_\alpha$ and project to $\mu$. For every $\varepsilon>0$, choose a compact $K\subset X$ with $\mu(K)>1-\varepsilon$. Then $\pi^{-1}K$ is compact and
\[
 \lambda_T(\pi^{-1}K)=\mu(K)>1-\varepsilon
 \qquad\text{for all }T>0.
\]
Thus $(\lambda_T)$ is uniformly tight. Hence any weak-$*$ limit of $(\lambda_T)$ as $T\to\infty$ is an $A$-invariant probability measure supported on $\widetilde\Omega_\alpha$ that projects to $\mu$.
\end{proof} 
We now claim that
\begin{equation}
\label{eq:omega_transport}
    b^k\widetilde\Omega_\alpha = \widetilde\Omega_{p^k\alpha}.
\end{equation}
Indeed, we note the identity
\begin{equation*}
 \begin{pmatrix}1&-\alpha\\0&1\end{pmatrix}
 \begin{pmatrix}1&0\\0&n\end{pmatrix}
 =\begin{pmatrix}1&0\\0&n\end{pmatrix}
 \begin{pmatrix}1&-n\alpha\\0&1\end{pmatrix}.
\end{equation*}
Hence
\begin{equation}
\label{eq:time_shift}
     b^ky_\alpha=(u^+(-\alpha),\diag(p^k,1))\Gamma_p
 =a_{-c_k}y_{p^k\alpha},
\end{equation}
where $c_k = \frac{k}{2}\log p$. The element $b$ commutes with $A$, and the fixed time shift $a_{-c_k}$ leaves the omega-limit set unchanged, whence the claim. In particular, if a measure $\sigma$ is supported on $\widetilde\Omega_\alpha$, then $(b^k)_*\sigma$ is supported on $\widetilde\Omega_{p^k\alpha}$.

\begin{proof}[Proof of Theorem~\ref{theo:cf_patterns}]
Choose $\sigma$ as in Lemma~\ref{lem:supported_lift}. By Theorem~\ref{cor:haar_cesaro}, there is a set $D$ of density one such that $\nu_k:=\pi_*((b^k)_*\sigma)\to m_X$ as $k\to\infty$ in $D$. Equations~\eqref{eq:omega_projection} and~\eqref{eq:omega_transport} give
\begin{equation*}
 \operatorname{supp}\nu_k\subseteq\Omega_{p^k\alpha}.
\end{equation*}
The continued fraction coding of the geodesic flow gives for every pattern $w$ a non-empty relatively compact open set $V_w\subset X$ such that
\begin{equation*}
\label{eq:word_box}
 \Omega_\beta\cap V_w\ne\varnothing
 \quad\Longrightarrow\quad
 w\text{ occurs infinitely often in the continued fraction expansion of }\beta.
\end{equation*}
See \cite{Series1985} and \cite[Appendix~A]{EFS2011} for details of the coding. 

Fix a pattern $w$. Since $V_w$ is non-empty and open,
\[
 \liminf_{\substack{k\to\infty\\k\in D}}\nu_k(V_w)\geq m_X(V_w)>0.
\]
Thus $\Omega_{p^k\alpha}$ meets $V_w$ for all sufficiently large $k\in D$, as required.
\end{proof}

For completeness, we include the standard deduction of Corollary~\ref{cor:littlewood} from Theorem~\ref{theo:cf_patterns}, applied to one-digit patterns.

\begin{proof}[Proof of Corollary~\ref{cor:littlewood}]
Fix an integer $m\geq1$. Theorem~\ref{theo:cf_patterns} gives $k$ such that the digit $m$ occurs infinitely often in the continued fraction expansion of $p^k\alpha$. A convergent immediately preceding such a digit has denominator $q$ with $q\|q\, p^k\alpha\|<1/m$. Taking $Q=p^kq$ gives
\[
 Q|Q|_p\|Q\alpha\|
 =q|q|_p\|q\, p^k\alpha \|<\frac{1}{m}.
\]
Letting $m\to\infty$ proves the corollary.
\end{proof}

\subsection{Statistics of the patterns}

We indicate how Corollary~\ref{cor:gauss-kuzmin} follows from Theorem~\ref{theo:padic_equidistribution}. Choose a $\mu$-generic point $x\in \Omega_\alpha$ and lengths $T_j \to \infty$. This means that $\frac{1}{T_j}\int_0^{T_j} \delta_{a_tx}\, dt \to \mu$ as $j \to \infty$. For each $j$, one can pick some time $s_j> j$ such that $a_{s_j +t}x_\alpha $ is $\frac{1}{j}$-close to $a_{t}x$ on the entire interval $0\leq t \leq T_j$. In particular, it follows that
\begin{equation*}
    \mu_j := \frac{1}{T_j}\int_{s_j}^{s_j+T_j} \delta_{a_tx_\alpha}\, dt \to \mu
\end{equation*}
as $j \to \infty$. We can similarly consider the lifted orbit measures
\begin{equation*}
    \sigma_j := \frac{1}{T_j}\int_{s_j}^{s_j+T_j} \delta_{a_ty_\alpha}\, dt
\end{equation*}
where $\pi_*\sigma_j = \mu_j$, and since $\pi$ is proper we can pass to a convergent subsequence, say with $\sigma_j \to \sigma$. The measure $\sigma$ is $A$-invariant, $\pi_*\sigma = \mu$, and $\supp \sigma \subset \widetilde\Omega_\alpha$. 

Recall from \eqref{eq:time_shift} that
\begin{equation*}
    b^ky_\alpha = a_{-c_k}y_{p^k\alpha},
\end{equation*}
where $c_k = \frac{k}{2}\log p$. Since $b$ commutes with $A$, for every fixed $k$ we obtain that
\begin{equation*}
    \pi_*((b_k)_*\sigma_j) = \frac{1}{T_j}\int_{s_j-c_k}^{s_j-c_k+T_j}\delta_{a_tx_{p^k\alpha}}\, dt \to \pi_*((b_k)_*\sigma).
\end{equation*}
By Theorem~\ref{theo:padic_equidistribution} we have a natural density one set $D$ such that $\pi_*((b_k)_*\sigma) \to m_X$ along $D$. By a diagonal argument, ensuring that we pick $j$ large enough so that the start time $s_j -c_k >0$, Corollary~\ref{cor:gauss-kuzmin} follows.

	\begingroup
    \hyphenpenalty=50
    \exhyphenpenalty=50
    \printbibliography[title=References]
    \endgroup

\end{document}